\documentclass[10pt,letterpaper,conference]{ieeeconf}

\usepackage{calrsfs,times}
\usepackage{graphicx,psfrag,enumerate}
 \usepackage{pstool}
\usepackage{amssymb}
\usepackage{amsmath}
\usepackage{calrsfs}
\usepackage{xcolor, epsfig, epstopdf, wrapfig}
\usepackage{paralist}
\usepackage{float,cuted}
\usepackage[colorlinks,bookmarksopen,bookmarksnumbered,citecolor=blue]{hyperref}
\usepackage{xargs}

\usepackage{todonotes}
\newcommandx{\iman}[2][1=]{\todo[linecolor=blue,backgroundcolor=blue!25,bordercolor=blue,size=small,author=Iman,#1]{#2}}

\newtheorem{definition}{Definition}
\newtheorem{theorem}{Theorem}
\newtheorem{remark}{Remark}
\newtheorem{example}{Example}
\newtheorem{lemma}{Lemma}

\newtheorem{assumption}{Assumption}

\usepackage{amsmath}
\usepackage{accents}
\newlength{\dhatheight}

\usepackage{color}
\usepackage[normalem]{ulem}
\usepackage{soul} 
\soulregister\cite7
\soulregister\ref7
\soulregister\pageref7

\IEEEoverridecommandlockouts
\allowdisplaybreaks

\begin{document}

\title{\LARGE \bf Global convergence and { asymptotic} optimality of the heavy ball
  method for { a class of} non-convex optimization { problems}\thanks{This research was supported in
    part under Australian Research Council's Discovery Projects funding
    scheme (DP190102158, DP200102945 and DP210102454). Part of this work was carried out
    during the first author's visit to the Australian National University.}}

\author{V. Ugrinovskii\thanks{V. Ugrinovskii is with the School of
    Engineering and IT, UNSW Canberra, Canberra, Australia, (e-mail:
  v.ugrinovskii@gmail.com).}
      \and I. R. Petersen \and I. Shames\thanks{I. R. Petersen and I. Shames
        are with CIICADA Lab, School of Engineering, The Australian National University, Canberra, ACT
        2601, Australia (e-mail: i.r.petersen@gmail.com;
        iman.shames@anu.edu.au).}}

\maketitle

\begin{abstract}
  In this letter we revisit the famous heavy ball method and study its
  global convergence { for a class of} non-convex problems { with sector-bounded
  gradient}. We characterize the parameters that render the method globally
convergent and yield the best $R$-convergence factor. We show that for { this}
family of functions, this convergence factor is superior to the factor
obtained from the triple momentum method.
\end{abstract}

\section{Introduction}
%
%
We consider the optimization problem
\begin{equation}
  \label{eq:66}
  \min_{x\in \mathbf{R}^n} f(x)
\end{equation}
 along with the well-known heavy ball (HB) method
\begin{eqnarray}
  \label{eq:HB}
  x_{t+1}=x_t-\alpha\nabla f(x_t)+\beta (x_t-x_{t-1}), \quad { t=0,1,\ldots,}
\end{eqnarray}
where $f:\mathbf{R}^n\to \mathbf{R}$ is differentiable and potentially non-convex and $\nabla f(x)$ denotes the gradient of $f(x)$. The parameters $\alpha$ and $\beta$ are positive scalar constants. Convergence
 properties of such iterative processes  have been widely studied; e.g.,
 see~\cite{OR-2000,Polyak-1987,Nesterov-2018}. It is also recognized that iterative processes of the form~(\ref{eq:HB}) can be
analyzed using tools from robust control
theory~\cite{LRP-2016,MSE-2020,SE-2021,VFL-2017}, since
equation~(\ref{eq:HB}) can be written in the form of a nonlinear dynamic
system of Lur\'e type.

The original local convergence properties of \eqref{eq:HB} as well as optimal choices for $\alpha$ and $\beta$ that yield the fastest $R$-convergence properties, see \cite{OR-2000}, were presented in \cite{Polyak-1987}. We refer to these optimal parameters as \emph{Polyak's parameters} for the rest of this paper.
Since then many researchers have contemplated the question of the global convergence of \eqref{eq:HB}.
For example, as stated in \cite{LRP-2016} for a certain family of strongly convex cost functions with Lipschitz gradient, the system  \eqref{eq:HB} with Polyak's parameters may not be  globally convergent\footnote{We will clarify this and revisit an example given in \cite{LRP-2016} later in this letter.}. In \cite{GFJ-2015}, intervals for the algorithm parameters are introduced that render \eqref{eq:HB} globally convergent. However,  that work does not focus on the relationship between the method's parameters and its convergence rate. More recently, a globally convergent first order optimization algorithm, \emph{the triple momentum method (TMM)}, for solving optimization problems with strongly convex cost function whose gradient is Lipschitz is introduced in \cite{VFL-2017}. It is further shown that one can choose the parameters for this algorithm to obtain the best hitherto convergence property among all globally convergent methods.

In this paper, we consider a more general class of functions compared to the aforementioned results. Specifically, inspired by the observations in \cite{NNG-2019}, we consider optimization problems with unimodal objective functions which may not be globally convex, and are only convex in the infinitesimally small region in the vicinity of
the minimum $x^*$. Specifically, we consider the problems whose cost functions satisfy the following assumption.

\begin{figure}[t]
  \centering
  \includegraphics[width=0.8\columnwidth]{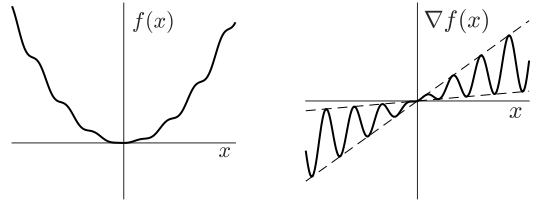}
  \caption{A non-convex objective function $f(x)\in\mathcal{F}_{m,L}^1$ where $f(x)=\frac{L-m}{4}\left( \frac{L +m}{L-m} x^2+\frac{2}{\omega ^2}\sin
    (\omega x)-\frac{2}{\omega }x \cos ( \omega x)\right)$,
    with positive scalars $L$, $m$, and $\omega$ such that $L>m$.}  \label{fig:ex1}
\end{figure}

\begin{assumption}\label{A.sector}
There exists a point $x^*$ such that the minimum of $f(x)$ is attained at $x^*$. Moreover, there exist constants $m$, $L$, $(L\ge m>0)$, such that for every
$x\in\mathbf{R}^n$,
\begin{equation}
  \label{eq:3}
  \left(m(x-x^*)-\nabla f(x)\right)^T\left(L(x-x^*)-\nabla f(x)\right)\le 0.
\end{equation}
The class of objective functions with this property will be denoted
$\mathcal{F}_{m,L}^1$.
\end{assumption}
An example of a function in $\mathcal{F}_{m,L}^1$ is depicted in Fig.~\ref{fig:ex1}.

We specifically employ the circle criterion for absolute stability of
discrete time Lur\'e systems with sector-bounded
nonlinearities to derive expressions for the parameters of the heavy ball method \eqref{eq:HB} that guarantee its
global convergence and optimize its $R$-convergence factor. To distinguish this choice of parameters from Polyak's parameters \cite{Polyak-1987}, we name them \emph{generalized heavy ball (GHB) parameters}.
Moreover, we investigate the global convergence of the triple momentum
method for non-convex functions using the
circle criterion. Specifically, for a given Lipschitz constant $L$, we
derive a bound on
the size of the sector, in terms of the ratio $L/m$, such
that TMM can be guaranteed to be globally converging for all
functions of class $\mathcal{F}_{m,L}^1$. We also
identify the cases where \eqref{eq:HB} with GHB parameters enjoys a better
$R$-convergence property compared to TMM.

\section{Background}
Before continuing further we briefly discuss the family of functions that satisfy Assumption~\ref{A.sector}.
Condition (\ref{eq:3}) is known as the
sector condition in the absolute stability theory~\cite{Khalil}. This
terminology owes to the graphical interpretation of
the inequality (\ref{eq:3}) in the case $x\in\mathbf{R}$ where~(\ref{eq:3}) implies that the graph of the derivative $f'(x)$ is
enclosed in the sector formed by the lines $m(x-x^*)$ and $L(x-x^*)$; see
Fig.~\ref{fig:ex1}.
It is worth to compare Assumption~\ref{A.sector}  with other assumptions on
the objective function $f(x)$ that appear in the literature. First,
since $f(x)$ attains its minimum at $x^*$, then $\nabla f(x^*)=0$.
When the function $f$ is two times differentiable, this allows us to
conclude from
(\ref{eq:3}) that the Hessian of $f$ at $x^*$, denoted
$\Delta=\nabla^2f(x^*)$, satisfies 
the often-quoted condition (e.g.,
see~\cite[Theorem~1, p.65]{Polyak-1987})
\begin{equation}
  \label{eq:52}
mI\le\Delta\le LI.
\end{equation}
We denote the class of such objective functions
$\mathcal{F}_{m,L}^2$.

The second connection is with the class $\mathcal{S}_{m_0,L_0}^{1,1}$ of
continuously differentiable strongly convex functions $f:\mathbf{R}^n\to
\mathbf{R}$ with Lipschitz continuous
gradient~\cite{LRP-2016,SE-2021,MSE-2020} such that
for all $x,y\in \mathbf{R}^n$
\begin{equation}
  \label{eq:10}
  m_0\|x-y\|^2\le (\nabla f(x)-\nabla f(y))^T(x-y)\le L_0\|x-y\|^2.
\end{equation}
It is shown in~\cite[Proposition~5]{LRP-2016} that any $f\in
\mathcal{S}_{m_0,L_0}^{1,1}$ satisfies the sector condition (\ref{eq:3})
with $m=m_0$, $L=L_0$, i.e., $\mathcal{S}_{m_0,L_0}^{1,1}\subset
\mathcal{F}_{m,L}^1$ when $m=m_0$, $L=L_0$.
 In this case,
$m$ characterizes the coercivity property of
the convex function $f$ while $L$ is the Lipschitz constant of
$\nabla f$. Of course,
Assumption~\ref{A.sector} describes a larger class of objective functions
than the class of convex functions captured by (\ref{eq:10}). It includes
functions which are not convex (see Fig.~\ref{fig:ex1}).
{  It follows from (\ref{eq:3}) that 
  every
  sector-bounded $f$ is an invex function~\cite{BM-1986}. When
$x\in\mathbf{R}$, every such function is pseudoconvex.}

For a strongly convex function $f\in\mathcal{S}_{m,L}^{1,1}$, $\kappa=L/m$
is known as the condition number, but for functions
$f\in\mathcal{F}_{m,L}^1$, it characterizes the relative `width' of the
sector containing $\nabla f$.
This ratio will help us to
compare our results with those obtained in the literature for strongly
convex functions.
%
%
%

\subsection{Heavy ball method as Lur\'e feedback control
  systems}
As pointed out in~\cite{LRP-2016,MSE-2020,SE-2021}, the iterative
process~(\ref{eq:HB}) can be written in the form of a
nonlinear dynamic system
  \begin{eqnarray}
    \label{eq:51}
&&X_{t+1}=AX_t+BU_t, \\
&& Y_t=CX_t,\quad U_t=-\phi(Y_t), \nonumber
  \end{eqnarray}
whose state, output and nonlinearity are respectively
$X_t=\begin{bmatrix}x_{t-1}^T &  x_t^T
\end{bmatrix}^T$, $\phi(Y_t)=   \nabla f(Y_t).$
The matrices $A\in \mathbf{R}^{2n\times 2n}$, $B\in   \mathbf{R}^{2n\times n}$, and $C\in \mathbf{R}^{n\times 2n}$ are defined as $A=A_0\otimes I_n, \quad B=B_0\otimes I_n, \quad  C=C_0\otimes I_n$ with
\begin{align*}
  A_0 = \begin{bmatrix}
    0 & \quad 1\\ -\beta & \quad 1+\beta
\end{bmatrix}, \quad
B_0= \begin{bmatrix}
  0 \\ \alpha
\end{bmatrix}, \quad C_0= \begin{bmatrix}
  0 & \quad 1
\end{bmatrix}.
\end{align*}
Similarly, $X^*=\begin{bmatrix}1 &~ 1
\end{bmatrix}^T\otimes x^*$ is the unique
equilibrium of~(\ref{eq:51}). { We refer the reader
  to~\cite{LRP-2016} for a detailed introduction to analysis and design
  of iterative optimization algorithms using tools from robust control theory.}

The system~(\ref{eq:51}) is the feedback control system of Lur\'e
type~\cite{HC-2008,Khalil}. There are a number of results in the
literature concerned with stability of Lur\'e systems.
In this note we will
make use of the multivariate version of the so-called Tsypkin's circle
criterion~\cite{HC-2008} to obtain sufficient conditions for global
asymptotic convergence of the Lur\'e systems of the form~(\ref{eq:51})
subject to the sector condition~(\ref{eq:3}); { also
  see~\cite{HB-1994}}.
Admittedly, the circle criterion is the simplest and the least accurate in terms of capturing properties of the function $f$. However, it proves to be sufficient for the purpose of this letter which is to demonstrate the applicability of the heavy ball method beyond the class of strongly convex objectives.

%

\subsection{$R$-convergence}
Here we recall some basic definitions which formalize the notion of the rate of
convergence of an iterative process.

\begin{definition}[\cite{OR-2000},p.288]\label{Def.R-convergence}
  Let $\{x_t\}$ be a sequence that converges to $x^*$. Then the number
$    r_{\{x_t\}}=\limsup_{t\to \infty} \|x_t-x^*\|^{1/t}
$ 
is the \emph{root-convergence factor}, or \emph{$R$-factor} of $\{x_t\}$.
If $\mathcal{I}$ is an iterative process with limit point
$x^*$, and $\mathfrak{C}(\mathcal{I},x^*)$ is the set of all sequences generated by
$\mathcal{I}$ which converge to $x^*$, then
$  r_{\mathcal{I}}=\sup\{r_{\{x_t\}}: \{x_t\}\in \mathfrak{C}(\mathcal{I},x^*)\}$
is the $R$-factor of $\mathcal{I}$ at $x^*$.
\end{definition}

{ \begin{remark} The
  $R$-factor characterizes an asymptotic convergence rate of $\{x_t\}$.
\end{remark}}

\begin{definition}\label{Def.convergence}
Let $x_t$ denote the vector of iterates of the algorithm (\ref{eq:HB})
at step $t$ initiated at time $0$ with $x_0,x_{-1}$.
The algorithm is said to converge globally asymptotically to $x^*$
if its iterates have the following properties:
\begin{enumerate}[(a)]
\item
for every $\epsilon>0$ there exists
$\delta=\delta(\epsilon)>0$ such that $\|x_{-1}-x^*\|<\delta$,
$\|x_0-x^*\|<\delta$ imply that
$\|x_t-x^*\|<\epsilon$ for all $t\ge 0$; and
\item
    $\lim_{t\to \infty}\|x_t-x^*\|=0$ uniformly in $x_0$, $x_{-1}$
    for all $x_0,x_{-1}\in \mathbf{R}^n$.
\end{enumerate}
\end{definition}

\section{Global convergence analysis}
In this section we analyze two optimization
methods~(\ref{eq:HB}), the heavy ball method due to
Polyak~\cite{Polyak-1964a,Polyak-1987}
and the triple momentum method of~\cite{VFL-2017} when applied to non-convex
objective functions from class  $\mathcal{F}_{m,L}^1$ and are initiated
arbitrarily far from $x^*$. Both methods are
two-step first order `accelerated gradient' methods, manifesting excellent
convergence properties. The heavy ball method with Polyak's parameters achieves the
smallest $R$-factor  when it is initiated
sufficiently close to $x^*$ and $f \in \mathcal{F}_{m,L}^2$~\cite{Polyak-1987}. Hence,  it
is the `asymptotically fastest' method.

Our analysis employs the circle criterion for absolute stability of
{ general}
discrete-time { time-varying} Lur\'e systems with sector-bounded
nonlinearity given in~\cite[p.841]{HC-2008}. We present it in the form
{ specialized for}
the setting of the algorithm~(\ref{eq:HB}) and
Lur\'e system~(\ref{eq:51}) with a nonzero point of attraction.
Let $Q(z)$ be the transfer function of the linear part of the Lur\'e
system~(\ref{eq:51}). Then $Q(z)=Q_0(z)\otimes I_n$, where
\begin{equation}
  \label{eq:49}
  Q_0(z)\sim
  \left[
    \begin{array}{c|c}
      A_0 & B_0 \\ \hline C_0 & 0
    \end{array}
  \right].
\end{equation}

\begin{lemma}[Circle criterion, \cite{HC-2008}, p.841]\label{L-absstab}
Suppose that the transfer function
\begin{equation}
  \label{eq:4}
H(z)=[1+ L Q_0(z)][1+mQ_0(z)]^{-1}\otimes I_n
\end{equation}
is strict positive real. Then for any function $f\in \mathcal{F}_{m,L}^1$,
the algorithm (\ref{eq:HB}) globally asymptotically converges to $x^*$.
\end{lemma}

\subsection{Global convergence of \eqref{eq:HB} with Polyak's parameters}

It is well
known~\cite{Polyak-1964a,Polyak-1987} that
when $f\in\mathcal{F}_{m,L}^2$ and
\begin{eqnarray}
  \label{eq:14}
  0\le \beta < 1, \quad 0<\alpha<2(1+\beta) L^{-1},
\end{eqnarray}
one can find an $\epsilon>0$ such that for
any $x_0$, $x_{-1}$ such that $\|x_{-1}-x^*\|\le \epsilon$, $\|x_0-x^*\|\le
\epsilon$, the sequence of
iterates $x_t$ generated by (\ref{eq:HB}) converges to $x^*$ exponentially fast, $  \|x_t-x^*\|\le c(\delta)(\rho+\delta)^t$,
where  $0\le \rho<1$, $0<\delta<1-\rho$.
 Under these conditions, the minimal $\rho$ is achieved
when
\begin{equation}
  \label{eq:2}
        \alpha= \alpha_{\mathrm{P}} \triangleq \frac{4}{(\sqrt{L}+\sqrt{m})^2}, \quad
        \beta=\beta_{\mathrm{P}}\triangleq \frac{(\sqrt{L}-\sqrt{m})^2}{(\sqrt{L}+\sqrt{m})^2}
\end{equation}
and is equal to
\begin{align}
  r_{\mathrm{P}} \triangleq \frac{\sqrt{L}-\sqrt{m}}{\sqrt{L}+\sqrt{m}}.
\end{align}
The parameters $\alpha_{\mathrm{P}}$ and $\beta_{\mathrm{P}}$ are the  Polyak's parameters discussed in the introduction.

According to~\cite{LRP-2016}, for a strongly convex function $f$
subject to the sector condition~(\ref{eq:3}), convergence of the heavy
ball algorithm (\ref{eq:HB}), (\ref{eq:2}) can be guaranteed when $\kappa=L/m$ is
approximately equal to $6$ or less. This threshold is calculated numerically \cite{LRP-2016} and represents a point where the LMIs used to characterize the global stability of (\ref{eq:HB}) fail to be satisfied. The next theorem gives a precise meaning to this observation. It derives a range of ratios $L/m$ ($L$, $m$
are the sector bound constants from~(\ref{eq:3})) for which
the heavy ball algorithm (\ref{eq:HB}) with Polyak's parameters converges globally for all functions of
class $\mathcal{F}_{m,L}^1$ in the sense of Definition~\ref{Def.convergence}.

\begin{theorem}
  \label{T.circle-criterion}
If
\begin{equation}
  \label{eq:6}
 L/m< \kappa_0\triangleq 3+2\sqrt{2},
\end{equation}
then the heavy ball method (\ref{eq:HB}), (\ref{eq:2}) globally
asymptotically converges to $x^*$ for any $f\in\mathcal{F}_{m,L}^1$.
\end{theorem}
\begin{figure}[t]
  \centering
  \includegraphics[width=0.75\columnwidth]{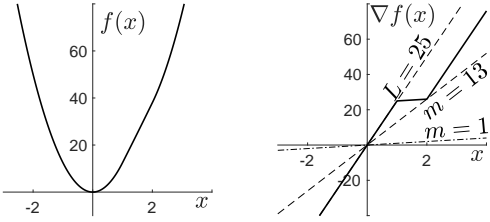}
  \caption{The function $f(x)$ in Example~\ref{Example-Lessard}.}
  \label{fig:ex2}
\end{figure}
\begin{remark}
  The threshold $\kappa_0=3+2\sqrt{2}$ on the ratio $L/m$ stated in
  Theorem~\ref{T.circle-criterion}  is
  close to the approximate threshold of 6 obtained in~\cite{LRP-2016} using
  numerical search. This is not surprising as the circle criterion is the approach that underpins the stability analysis of \cite{LRP-2016} in deriving this threshold. \hfill$\blacksquare$
\end{remark}
\emph{Proof of Theorem~\ref{T.circle-criterion}: }
With $\alpha$ and $\beta$ given in (\ref{eq:2}), the transfer function
$H(z)$ in equation (\ref{eq:4}) becomes $
H(z)=\left(\frac{z+\beta^{1/2}}{z-\beta^{1/2}}\right)^2 I_n$.
According to Lemma~\ref{L-absstab},
to establish the claim of the theorem it suffices to show that if $L/m
<\kappa_0=3+2\sqrt{2}$, then $H(z)$ is strictly positive
 real. This requires us to show that
$
 \left|\arg \left(\frac{e^{j\omega}+\beta^{1/2}}{e^{j\omega}-\beta^{1/2}}\right)\right|<\frac{\pi}{4}$ for all
$\omega \in[-\pi,\pi]$.
Equivalently,
it suffices to show that
\begin{eqnarray}
  \label{eq:9}
\min_{\omega \in[-\pi,\pi]} \cos\left(\arg
   \left(\frac{e^{j\omega}+\beta^{1/2}}{e^{j\omega}-\beta^{1/2}}\right)\right)
>\frac{1}{\sqrt{2}}.
\end{eqnarray}
However, 
$
\min_{\omega \in[-\pi,\pi]} \cos\left(\arg
   \left(\frac{e^{j\omega}+\beta^{1/2}}{e^{j\omega}-\beta^{1/2}}\right)\right)=
 \frac{1-\beta}{1+\beta}$.
Using the definition of $\beta$ in~(\ref{eq:2}), it is straightforward to verify
that $\frac{1-\beta}{1+\beta}>\frac{1}{\sqrt{2}}$ and~$H(z)$ is strictly
positive real if and only if $L/m<\kappa_0$. Thus, the conditions of the
circle criterion in Lemma~\ref{L-absstab} are satisfied under
condition~(\ref{eq:6}), and the  statement of the theorem follows from
Lemma~\ref{L-absstab}.
\hfill$\Box$
\begin{figure}[t]
  \centering
  \includegraphics[width=1\columnwidth]{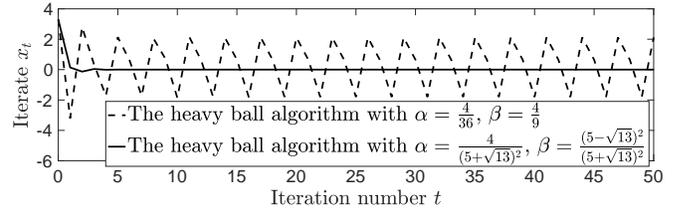}
  \caption{The trajectories of the two versions of the heavy ball algorithm
    for the function in
    Example~\ref{Example-Lessard}.}
  \label{fig:ex2-trajectories}
\end{figure}
\begin{example}
  \label{Example-Lessard}
Reference~\cite{LRP-2016} gives an example of a strongly convex function
$f(x)$ of one-dimensional variable $x$ that leads to a non-convergent heavy
ball method (see Fig.~\ref{fig:ex2}):
\begin{eqnarray}
  \label{eq:11}
&&  f(x)=
  \begin{cases}
    \frac{25}{2}x^2, & x<1, \\
    \frac{1}{2}x^2+24x-12, & 1\le x<2, \\
    \frac{25}{2}x^2-24x+36, & x\ge 2,
  \end{cases} 
\end{eqnarray}
The function has minimum at $x^*=0$. It satisfies (\ref{eq:10}) with
$m_0=1$, $L_0=25$, which implies that it also satisfies the sector
condition (\ref{eq:3}) with $m=1$, $L=25$. It was observed
in~\cite{LRP-2016} that using these
values of $m$ and $L$ in (\ref{eq:2}) to define the Polyak's parameters leads to an algorithm which,
when initiated within a certain range of initial conditions,
converges to a limit cycle, instead of the minimum
point $x^*=0$. The values of $\alpha$ and $\beta$ for which this behavior
was observed were $\alpha=\frac{4}{36}$,
$\beta=\frac{4}{9}$.

However, from Fig.~\ref{fig:ex2} one can see readily that the function
$f(x)$ in this example also satisfies the sector condition (\ref{eq:3}) with
$m=13$, $L=25$. These values of $m$, $L$ render the condition~(\ref{eq:6})
of Theorem~\ref{T.circle-criterion} true, thus the heavy ball
algorithm~(\ref{eq:HB}) with Polyak's parameters
$\alpha=\frac{4}{(5+\sqrt{13})^2}$,
$\beta=\frac{(5-\sqrt{13})^2}{(5+\sqrt{13})^2}$
obtained using $m=13$, $L=25$ must globally asymptotically
converge to $x^*=0$.
Fig.~\ref{fig:ex2-trajectories} confirms this.
Both trajectories in this figure are initiated with the same initial condition
$x_{-1}=x_0=3.3$.

This example highlights the importance of deciding what values $m$ and $L$ are to be used for tuning the heavy ball method. The values of $m$ and $L$ such that the condition~(\ref{eq:10}) is satisfied may be conservative for tuning the heavy ball algorithm, even for
convex functions. However, using values of $m$ and $L$ such that the sector condition~(\ref{eq:3}) is satisfied,
may be less conservative. This is not surprising
since these conditions are not equivalent: (\ref{eq:10})
implies~(\ref{eq:3}) but not vice versa. Thus, one should exercise caution when quoting this example as a counter-example for global convergence of \eqref{eq:HB}. 
\end{example}
\subsection{Global convergence of TMM}
The triple momentum method \cite{VFL-2017} is given by
\begin{align}
  \label{eq:72}
  x_{t+1}&=(1+\beta)x_t-\beta x_{t-1}-\alpha \nabla f(y_t), \nonumber\\
       y_t&=(1+\gamma)x_t-\gamma x_{t-1}, \nonumber \\
       \eta_t&=(1+\delta)x_t-\delta x_{t-1}.
\end{align}
The coefficients $\alpha$, $\beta$, $\gamma$ and $\delta$ were
defined in \cite{VFL-2017} as
\begin{equation}
  \label{eq:64}
(\alpha,\beta,\gamma,\delta)=\left(\frac{1+\rho}{L},
\frac{\rho^2}{2-\rho},\frac{\rho^2}{(1+\rho)(2-\rho)},\frac{\rho^2}{1-\rho^2}
  \right),
\end{equation}
where $\rho=1-\sqrt{\kappa^{-1}}$, $\kappa$ is the condition number of
a strongly convex function, and $L$ is the Lipschitz constant of its gradient;
see~(\ref{eq:10}). As noted, any such function belongs to the class
$\mathcal{F}_{m,L}^1$ with $m=L/\kappa$.
The process~(\ref{eq:72}) can be
represented as a Lur\'e system of the
form~(\ref{eq:51}) with the state $X_t=[x_{t-1}^T~x_t^T]^T$, the output
$Y_t=[-\gamma I_n~(1+\gamma)~I_n]X_t$ and the additional performance output
$\eta_t=[-\delta I_n~(1+\delta)~I_n]X_t$ which is not fed back
into~(\ref{eq:51}).
Therefore global asymptotic convergence of this
method can be analyzed using the circle criterion, and the assumption of
strong convexity of $f$ can be relaxed to $f\in\mathcal{F}_{m,L}^1$, where
$m=L/\kappa$, and $\kappa>0$ can be regarded as a given constant.

Consider the polynomial
$ 
  \chi(\rho)=8-\rho-8\rho^2-14\rho^3-\rho^5.
$ 
This polynomial is monotone decreasing in $\rho$ when $\rho>0$, and it has
the unique real root $\rho_0\approx 0.650307$. Define $\kappa_{\mathrm{TM}}\triangleq(1-\rho_0)^{-2}\approx 8.1776.$
\begin{theorem}
  \label{T.TM}
Let $\kappa$ be a positive constant, and $f$ be an arbitrary function in $\mathcal{F}_{m,L}^1$, where $m=L/\kappa$.
If $\kappa<\kappa_{\mathrm{TM}}$, then each sequence of iterates $x_t$,
$y_t$ and $\eta_t$ generated by~(\ref{eq:72}) globally asymptotically
converges to $x^*$.
\end{theorem}

\emph{Proof: }
With $\alpha$, $\beta$, $\gamma$ and $\delta$ defined in~(\ref{eq:64}), the
transfer function $H(z)$ defined in Lemma~\ref{L-absstab} is
$ H(z)=
\frac{z (\rho +z)}{(z-\rho ) \left(z-\rho ^2\right)}I.
$ 
We now show that $H(z)$ is strict positive real for
$\kappa<\kappa_{\mathrm{TM}}$. 
For this, it suffices to show that for all $\omega \in[-\pi,\pi]$,
\begin{eqnarray}
\psi(\omega) &\triangleq& \mathrm{Re}
\left[e^{j\omega}(\rho
  +e^{j\omega})(e^{-j\omega}-\rho)(e^{-j\omega}-\rho^2)\right] \nonumber \\
&=&\rho ^4 \cos \omega -2 \rho ^3 \sin ^2\omega -\rho ^2 \cos \omega
    -\rho ^2+1 >0. \qquad
\label{eq:75}
\end{eqnarray}

When $\rho\in (0,\sqrt{5}-2)$, the function $\psi(\omega)$ has minimum at
$\omega=0$, and $\psi(0)=(1-\rho^2)^2>0$ for any $\rho$ in this interval.

When $\rho\in[\sqrt{5}-2,1)$, the function $\psi(\omega)$ has two
minimums at $\omega=\pm\cos^{-1}\left(\frac{1-\rho^2}{4\rho}\right)$. The
function is even, therefore it has the same value of $\frac{1}{8}\chi(\rho)$ at those points.
Since $\chi(\rho)$ is monotone decreasing and $\chi(\rho)\le 0$ for
$\rho\ge \rho_0$, we conclude that $H(z)$ is strict
positive real in this case if and only if $\rho< \rho_0$.

Combining these two cases, we note that the conditions of Lemma~\ref{L-absstab}
are satisfied when $\kappa<\kappa_{\mathrm{TM}}$, and the claim of the theorem
follows from that lemma. \hfill$\Box$

Theorem~\ref{T.TM} shows that TMM tuned for a strongly convex
function $f\in\mathcal{S}_{m_0,L}^{1,1}$ with the condition number
$\kappa<\kappa_{\mathrm{TM}}$ remains globally converging for non-convex
functions of the class $\mathcal{F}_{m,L}^1$ with $m=m_0=L/\kappa$.

\section{Global convergence of the heavy ball method with GHB parameters and the corresponding $R$-convergence factor}
Example~\ref{Example-Lessard} motivates us to address the question raised
in~\cite{LRP-2016} as to whether the heavy ball method converges globally
over the class of functions captured by the sector condition~(\ref{eq:3})
for some choice of $\alpha$ and $\beta$ than other the Polyak's parameters
(\ref{eq:2}). 
Using the circle
criterion allows us to answer this question, without making
\emph{a priori} assumptions about the ratio $\kappa=L/m$.

Introduce the function
\begin{equation}
  \label{eq:12}
  \bar\alpha(\beta) = \begin{cases} 2(1+\beta)L^{-1}, & \beta \le
    \left(\sqrt{\frac{L}{m}}-\sqrt{\frac{L}{m}-1}\right)^2, \\
\frac{2(1-\beta)^2}{(1+\beta)
  (L+m)-4\sqrt{\beta Lm}}, & \beta > \left(\sqrt{\frac{L}{m}}-\sqrt{\frac{L}{m}-1}\right)^2;
\end{cases}
\end{equation}

\begin{theorem}\label{T.SPR}
For any $L$, $m$ such that $L>m>0$, if
\begin{equation}
  \label{eq:17}
  0\le \beta<1, \quad 0<\alpha<\bar\alpha(\beta),
\end{equation}
then the algorithm (\ref{eq:HB}) globally
asymptotically converges to $x^*$ for all functions $f\in \mathcal{F}_{m,L}^1$.
\end{theorem}




\begin{figure*}[!t]
\normalsize
\begin{equation}
  \label{eq:37}
\beta_0(\kappa)=
\frac{\kappa  \left(\kappa  \left(-\sqrt{2} \sqrt{\frac{(\kappa -1) \left(\left(\sqrt{\frac{\kappa -8}{\kappa }}+1\right) \kappa ^2+\left(7 \sqrt{\frac{\kappa -8}{\kappa }}-5\right) \kappa +12\right)}{\kappa ^3}}+\sqrt{\frac{\kappa -8}{\kappa }}+1\right)-\sqrt{\frac{\kappa -8}{\kappa }}+7\right)^2}{16 (\kappa +1)^2}
\end{equation}
\hrulefill
\vspace*{4pt}
\end{figure*}

\emph{Proof: }
According to Lemma~\ref{L-absstab}, it suffices to show that
\begin{equation}
  \label{eq:13}
  H(z)=\frac{z^2+\left(L\alpha -1-\beta\right)z+\beta}
{z^2+\left(m\alpha-1-\beta\right)z+\beta}I_n
\end{equation}
is strictly positive real when $\alpha$, $\beta$ satisfy (\ref{eq:17}). For
this, we use the criterion for strict positive
realness due to \v{S}iljak; see~\cite[Theorem~2]{Siljak-1973}.

Let $p(z)$, $q(z)$ denote the denominator and the numerator of $H(z)$.
Since $\bar\alpha(\beta)< 2(1+\beta)L^{-1}$ in the interval
$\beta\in\left(\left(\sqrt{L/m}-\sqrt{L/m-1}\right)^2,1\right]$,
one can use the Jury criterion to confirm that the polynomial
\begin{equation}
  \label{eq:19}
\frac{1}{2}(p(z)+q(z))=z^2+\left(\frac{m+L}{2}\alpha-1-\beta\right)z+\beta
\end{equation}
is Schur stable under conditions (\ref{eq:17}). This validates
condition (i) of Theorem~2 in~\cite{Siljak-1973}.

Condition (ii) of that theorem requires that the
polynomial
$
g(z)=\frac{1}{2}z^2(p(z)q(z^{-1})+p(z^{-1})q(z))$
must have exactly two roots inside and two roots outside the unit
circle. To validate this requirement, we note that $g(z)$ is a self-inversive polynomial\footnote{A polynomial $g(z)$ of degree $l$ is
self-inversive if $g(z)=z^{l}g(z^{-1})$.}, hence it can be written as
$g(z)=b_2z^4+b_1z^3+2b_0z^2+b_1z+b_2$. We now use the modified Schur-Cohn-Marden method for enumerating the zeros of a self-inversive polynomial
within the unit circle~\cite{BM-1952}. According to this method,
$g(z)$ has exactly two zeros inside and two zeros outside the unit disk if
and only if $
  b_1^2-\left(b_0+b_2\right)^2<0$
and either $
  16b_2^2-b_1^2<0$,
or
$
  16b_2^2-b_1^2>0$ and $
\left(16b_0b_2-3b_1^2\right)^2>\left(16b_2^2-b_1^2\right)^2.
$
The 
analysis of these conditions shows that they are equivalent to
condition (\ref{eq:17}).
\hfill$\Box$

\begin{remark}
  When $\beta\le \left(\sqrt{L/m}-\sqrt{L/m-1}\right)^2$,
  the upper bound on $\alpha$ in~(\ref{eq:17}) is the same as that given in
  Polyak's theorem on stability of the heavy ball method~\cite[Theorem~1,
  p.65]{Polyak-1987}; see~(\ref{eq:14}). Therefore, for this range of the
  inertia parameter $\beta$, Theorem~\ref{T.SPR} complements Polyak's
  theorem by showing that the
  heavy ball method converges globally. When $\beta>
  \left(\sqrt{L/m}-\sqrt{L/m-1}\right)^2$, Theorem~\ref{T.SPR} guarantees global
  convergence only when the stepsize parameter $\alpha$ is less than
  $\bar\alpha(\beta)$, i.e., in the smaller range of $\alpha$ than that
  established by Polyak's theorem. \hfill$\blacksquare$
\end{remark}
\begin{remark}
  In~\cite{GFJ-2015}, it is shown that for a strongly convex
  Lipschitz-continuous  $f$, if $
    \beta\in(0,1)$ and  $\alpha\in\left (0,\frac{2(1-\beta)}{L+m}\right ),
  $
  then \eqref{eq:HB} converges with an $R$-factor less than 1.
  We now compare this result with  Theorem~\ref{T.SPR}.
  It is easy to see that $\frac{2(1-\beta)}{L+m}<2(1+\beta)L^{-1}$. Also,
  when
  $\left(\sqrt{\frac{L}{m}}-\sqrt{\frac{L}{m}-1}\right)^2
  <\beta\le \frac{4Lm}{(L+m)^2}$, then $
    \bar\alpha(\beta)\ge \frac{2(1-\beta)}{L+m}$.
  When $\beta \in(\frac{4Lm}{(L+m)^2},1)$, the opposite inequality holds.
  These inequalities allow us to compare the conditions in \cite{GFJ-2015}  with the
  result of Theorem~\ref{T.SPR} and condition~(\ref{eq:17}) of that
  theorem. We conclude that for $\beta\in (0,\frac{4Lm}{(L+m)^2})$, our
  result provides a wider range of $\alpha$ under which the GHB method is
  globally convergent. This conclusion holds even when $f$ is not
  strongly convex. \hfill$\blacksquare$
\end{remark}
Next, we derive an upper bound on the optimal $R$-factor of the globally
converging algorithm (\ref{eq:HB}) for the parameter choice of Theorem~\ref{T.SPR}, which we denote $r_{\mathrm{GHB}}$.
This will be achieved by optimizing the $R$-factor of the
iterative process~(\ref{eq:HB}) over the parameter region
defined by the global convergence condition (\ref{eq:17}).
To derive such a bound, we restrict attention to
 objective functions $f\in\mathcal{F}_{m,L}^2$.

To present this result, some notation is needed. Let $\rho(A(\Delta))$ denote
the spectral radius of the matrix
\begin{equation}
  \label{eq:38}
A(\Delta)=  \left[
    \begin{array}{cc}
      0 & I \\
      -\beta I~ & ~(1+\beta)I- \alpha \Delta
    \end{array}
  \right], \quad \Delta=\nabla^2 f(x^*).
\end{equation}
Also, introduce $
  \gamma_m=\gamma_m(\alpha,\beta)=m\alpha-1-\beta$, and $
 \gamma_L=\gamma_L(\alpha,\beta)=L\alpha-1-\beta$.
It can be shown that for every $\kappa=L/m>8$ the function $
\eta(\beta)
\triangleq\frac{-\gamma_m(\bar\alpha(\beta),\beta)+\sqrt{\gamma_m(\bar\alpha(\beta),\beta)^2-4\beta}}{2}$
has a unique local minimum in the interval
$\beta\in ((\sqrt{\kappa}-\sqrt{\kappa-1})^2,1)$. This minimum is attained at
$\beta=\beta_0(\kappa)$ defined in equation (\ref{eq:37}) shown at the top
of the page.
Next, let us define $
  \nu(\kappa)\triangleq\frac{2-\sqrt{2\sqrt{\kappa}+3-\kappa}}{\sqrt{\kappa }-1}$.
This quantity is well defined in the interval $\kappa\in
[\kappa_0,9]$ where $\kappa_0$ is the constant defined in~(\ref{eq:6}).
Finally, it can be shown that in the interval $\kappa\in[8,9]$ the equation $
\eta(\beta_0(\kappa))=\nu(\kappa)$ has a unique solution $\bar\kappa \approx 8.2975$.

\begin{theorem}
  \label{T.stability-margin}
Let $\mathcal{A}_{\mathrm{GHB},m,L}$ denote the set of
all pairs $(\alpha, \beta)$ for which the algorithm~(\ref{eq:HB}) converges
globally for all $f\in\mathcal{F}_{m,L}^2$. Then $
\inf_{\mathcal{A}_{\mathrm{GHB},m,L}}\sup_{\mathcal{F}_{m,L}^2} r_{\mathrm{GHB}} \le r^*$,
where the constant $r^*$ is defined by
\begin{equation}
  \label{eq:22}
  r^*= \min \{\sup_\Delta\rho(A(\Delta)):
  \alpha\in(0,\bar\alpha(\beta)), \beta\in[0,1)\}.
\end{equation}
The supremum in~(\ref{eq:22}) is over the set of matrices
$\Delta$ which satisfy~(\ref{eq:52}). The constant $r^*$ has the following value
\begin{eqnarray}
\label{eq:30}
r^*=\begin{cases}
\frac{\sqrt{L}-\sqrt{m}}{\sqrt{L}+\sqrt{m}}, & 1<\frac{L}{m}\le \kappa_0, \\
\nu(L/m) & \kappa_0\le \frac{L}{m} < \bar\kappa, \\
\eta(\beta_0(L/m)),  & \frac{L}{m}\ge \bar\kappa.
\end{cases}
\end{eqnarray}
The minimum in~(\ref{eq:22}) is achieved when
\begin{eqnarray}
  \label{eq:41}
&& \alpha= \alpha^*\triangleq\begin{cases}
\frac{4}{(\sqrt{L}+\sqrt{m})^2}, & 1<\frac{L}{m}\le \kappa_0,\\
\bar\alpha\left(\nu^2(L/m)\right)
, & \kappa_0\le \frac{L}{m}< \bar \kappa, \\
\bar\alpha(\beta_0(L/m)), & \frac{L}{m}\ge \bar\kappa,
\end{cases}
\nonumber \\
&& \beta= \beta^*\triangleq\begin{cases}
\left(\frac{\sqrt{L}-\sqrt{m}}{\sqrt{L}+\sqrt{m}}\right)^2, & 1<\frac{L}{m}\le
\kappa_0, \\
\nu^2(L/m)
, & \kappa_0\le \frac{L}{m}< \bar\kappa, \\
\beta_0(L/m), & \frac{L}{m}\ge\bar\kappa.
\end{cases}
\end{eqnarray}
\end{theorem}

\emph{Proof: }
The starting point of the proof is the observation that the $R$-factor
of the iterative process~(\ref{eq:HB}) is determined by the
spectral radius of the matrix $A(\Delta)$ defined in~(\ref{eq:38});
see \cite[p.353]{OR-2000}.
Formally, we write this statement as $r_{\mathrm{GHB}} = \rho(A(\Delta))$.
Then, since $\{(\alpha, \beta):
\alpha\in(0,\bar\alpha(\beta)), \beta\in[0,1)\}\subseteq
\mathcal{A}_{\mathrm{GHB},m,L}$, $
\inf_{\mathcal{A}_{\mathrm{GHB},m,L}}\sup_{\mathcal{F}_{m,L}^2} r_{\mathrm{GHB}} \le r^*$ follows from this evaluation of  $r_{\mathrm{GHB}}$.

We now obtain a closed form expression for $r^*$, by
evaluating the minimum in~(\ref{eq:22}) and also obtain the parameters
$\alpha^*$, $\beta^*$ which attain that minimum. First, a closed form
expression for $\sup_{\Delta}\rho(A(\Delta))$ will be obtained.

Note that for any $\alpha$, $\beta$ which satisfy (\ref{eq:17})
the characteristic polynomial of the matrix $A(\Delta)$, i.e., $
   \det(z^2I+(\alpha\Delta-(1+\beta)I)z+\beta I)$,
is Schur stable (in fact, it is Schur stable in the larger
region defined by~(\ref{eq:14})). Equivalently, each polynomial
\begin{equation}
  \label{eq:20}
  z^2+(\alpha\lambda_i-1-\beta)z+\beta, \quad i=1,\ldots, n,
\end{equation}
is Schur stable; here $\lambda_i$ is the $i$-th eigenvalue of the matrix
$\Delta$, $\lambda_i\in[m,L]$.
Let $z_i^{\pm}$, $i=1, \ldots,n$, be the roots of the polynomials
(\ref{eq:20}). The scalar $\rho(A(\Delta))$ is the
radius of the smallest circle in the complex plane which contains all
$z_i^{\pm}$,
\begin{equation}
  \label{eq:32}
  \rho(A(\Delta))=\inf\{r\in[0,1): |z_i^\pm|\le r, i=1,\ldots,n \}.
\end{equation}
Using the Jury criterion, it is straightforward to obtain that $|z_i^\pm|\le
r$ if and only if
\begin{eqnarray}
\label{eq:34}
|\beta|\le r^2, \quad
r^2+\gamma_mr+\beta \ge 0, \quad
r^2-\gamma_Lr+\beta \ge 0. \quad
\end{eqnarray}
By computing the smallest $r\in[0,1)$ which satisfies (\ref{eq:34}) we
obtain that
\begin{eqnarray}
  \label{eq:26}
\lefteqn{\sup_\Delta\rho(A(\Delta))} && \nonumber \\
&=&\begin{cases}\sqrt{\beta}, & \text{\hspace{-5.5cm}when
      $\gamma_L^2<4\beta$, $\gamma_m^2<4\beta$ 
      ;} \\
\max\left(\sqrt{\beta},\frac{\gamma_L+(\gamma_L^2-4\beta)^{1/2}}{2}\right),
\\
&
\text{\hspace{-5.5cm}when
      $\gamma_L^2\ge 4\beta$, $\gamma_m^2<4\beta$ 
      ;} \\
\max\left(\sqrt{\beta},\frac{-\gamma_m+(\gamma_m^2-4\beta)^{1/2}}{2}\right),
\\
& \text{\hspace{-5.5cm}when
      $\gamma_L^2<4\beta$, $\gamma_m^2\ge 4\beta$ 
      ;} \\
\max\left(\sqrt{\beta},\frac{\gamma_L+(\gamma_L^2-4\beta)^{1/2}}{2},
\frac{-\gamma_m+(\gamma_m^2-4\beta)^{1/2}}{2}\right), \\
&  \text{\hspace{-5.5cm}when
      $\gamma_L^2\ge 4\beta$, $\gamma_m^2\ge 4\beta$ 
      .} 
\end{cases} \qquad \qquad 
\end{eqnarray}


To complete the proof, let us partition the region in the $\beta$,
$\alpha$ plane covered by condition~(\ref{eq:17}) into subregions,
according to the expressions for $\sup_\Delta\rho(A(\Delta))$ given
in~(\ref{eq:26}). This allows to minimize $\sup_\Delta\rho(A(\Delta))$ in each
subregion using the standard calculus tools, and so the optimal value over
the entire region can be obtained. Although these calculations are too
cumbersome to include here, they are straightforward and lead to the
expressions~(\ref{eq:30}),~(\ref{eq:41}).
\hfill$\Box$

As expected, the expressions for $\alpha$ and $\beta$ in~(\ref{eq:41})
and~(\ref{eq:30}) confirm the Polyak's optimal tuning of the heavy ball
algorithm when the ratio $L/m$ is within the interval
$(1,\kappa_0]$; this is consistent with Theorem~\ref{T.circle-criterion}. In this
case, the optimal pair $(\alpha^*,\beta^*)=(\alpha_{\mathrm{P}},\beta_{\mathrm{P}})$
lies in the interior of the region of global
convergence~(\ref{eq:17}).

\begin{figure}[t]
  \centering
  \includegraphics[width=1\columnwidth]{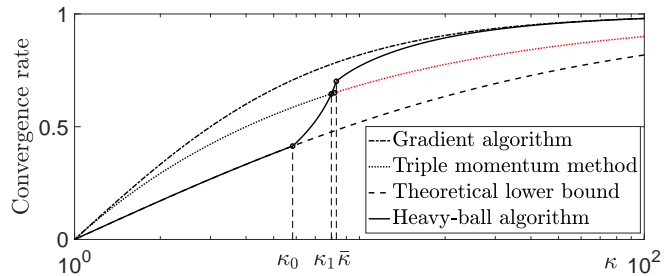}
  \caption{ Graphs of $r^*$ (the solid line) and the optimal convergence
    rates of the gradient algorithm (the dash-dot line) and the triple
    momentum method (the dotted line, extension beyond
    $\kappa_{\mathrm{TM}}$ is shown in red).
    The theoretical lower bound $r_{\mathrm{P}}$ is also shown (the dash line).}
  \label{fig:rstar}
\end{figure}

When $L/m\ge \kappa_0$, $(\alpha_{\mathrm{P}},\beta_{\mathrm{P}})$ defined by (\ref{eq:2}) lies
outside the region of global convergence~(\ref{eq:17}). The minimum of
$\sup_\Delta\rho(A(\Delta))$ over this region is
then attained on the boundary of the region, specifically on the line
determined by the second case of the expression in~(\ref{eq:12}). Nevertheless,
the heavy ball algorithm~(\ref{eq:HB}) with the GHB parameters selected according
to~(\ref{eq:41}) remains globally convergent and outperforms the optimally
tuned gradient algorithm, as Fig.~\ref{fig:rstar} shows. However, the benefits
of using the algorithm~(\ref{eq:HB}) with GHB parameters (\ref{eq:41}) over the gradient
algorithm diminish as the ratio $L/m$ increases.

Fig.~\ref{fig:rstar} also compares the obtained rate of global
convergence with the rate of convergence of the triple momentum
algorithm. It must be noted that the rate of global convergence of TMM for
functions in $\mathcal{F}_{m,L}^1$ is not known, however since  for any
$\kappa$, $\mathcal{S}_{m,L}^{1,1}\subset \mathcal{F}_{m,L}^1$ when $m=L/\kappa$,
the rate of convergence established for strongly convex functions
$\rho_{\mathrm{TM}}=1-1/\sqrt{\kappa^{-1}}$~\cite{VFL-2017} serves
is a lower bound on that rate. We conclude from Fig.~\ref{fig:rstar}
that in the interval $\kappa\in(1,\kappa_1]$, $\kappa_1=\frac{1}{2}\left(3
  \sqrt{7}+8\right)\approx 7.96863$ the triple momentum method converges
slower than \eqref{eq:HB} with GHB parameters.
Theorem~\ref{T.TM} does not guarantee that TMM converges globally for
functions in $\mathcal{F}_{m,L}^1$ with $m=L/\kappa$ when $\kappa\ge \kappa_{\mathrm{TM}}\approx 8.1776$. Therefore, our comparison is only valid
for $\kappa< \kappa_{\mathrm{TM}}$.



\section{Conclusions and Future Work}
We introduced the GHB parameters that lead to a globally convergent heavy
ball method while simultaneously minimizing an upper bound on the
$R$-convergence of the algorithm for solving a class of non-convex
optimization problems. We showed that the triple momentum method is not
guaranteed to be globally converging for all the functions in the same
class and identified a family of functions in the class that the heavy ball
method with GHB parameters enjoys a faster convergence property than the
triple momentum method. { Possible future directions include
  extending the result to constrained problems} and the cases where a
stochastic approximation of the gradient is available. { It is
  also interesting to find constructive conditions on $f$ which
  imply~(\ref{eq:3}).} Another avenue of research is to use more
sophisticated nonlinear stability analysis tools to obtain tighter  upper
bounds for the $R$-factor of the algorithm.


\newcommand{\noopsort}[1]{} \newcommand{\printfirst}[2]{#1}
  \newcommand{\singleletter}[1]{#1} \newcommand{\switchargs}[2]{#2#1}

\end{document}